\documentclass[10pt,reqno]{amsart}

\usepackage{amsmath,amsthm,amssymb,thmtools,tikz-cd,comment,hyperref,cleveref,enumitem,breqn,float,caption,subfig,ytableau,tikz,cases,kotex,multirow}
\usepackage{xcolor}
\usetikzlibrary{positioning, arrows.meta, calc}

\usetikzlibrary{decorations.pathreplacing,
	calligraphy,
	matrix}

\allowdisplaybreaks
\theoremstyle{plain}
\newtheorem{theorem}{Theorem}[section]
\newtheorem{lemma}[theorem]{Lemma}
\newtheorem{proposition}[theorem]{Proposition}

\newtheorem{thmalphabetintro}{Theorem}

\newtheorem{thmalphabetmaintext}{Theorem}

\theoremstyle{definition}
\newtheorem{definition}[theorem]{Definition}

\newtheorem{conjecture}[theorem]{Conjecture}

\newtheorem{open problem}[theorem]{Open Problem}
\newtheorem{remark and notation}[theorem]{Remark and Notation}
\newtheorem{remark and definition}[theorem]{Remark and Definition}
\newtheorem{definition and notation}[theorem]{Definition and Notation}
\newtheorem{notation and convention}[theorem]{Notation and convention}
\newtheorem{convention and notation}[theorem]{Convention and notation}

\def \p {\mathbb{P}}

\def \z {\mathbb{Z}}

\def \q {\mathbb{Q}}
\def \k {\mathbb{K}}

\def \o {\mathcal{O}}
\def \I {\mathcal{I}}
\def \f {\mathcal{F}}
\def \l {\mathcal{L}}
\def \C {\mathcal{C}}
\def \H {\mathcal{H}}

\def \sm {\textnormal{sm}}

\def \sing {\operatorname{sing}}
\def \leng {\operatorname{length}}
\def \supp {\operatorname{Supp}}

\def \bl {\operatorname{Bl}}

\def \codim {\operatorname{codim}}

\def \reg {\operatorname{reg}}

\def \rat {\operatorname{rat}}
\def \bs {\operatorname{Bs}}
\def \fac {\operatorname{fac}}

\def \cm {\operatorname{CM}}
\def \gen {\operatorname{gen}}

\def \la{\langle}
\def \ra{\rangle}

\def \vmd[#1]{\text{VMD$^{#1}$}}
\def \dpv[#1]{\text{dP$^{#1}$}}
\def \acm[#1,#2]{\text{ACMD$^{#1}_{#2}$}}
\def \propp[#1,#2]{\text{$P_{{#1},{#2}}$}}
\def \propa[#1,#2]{\text{$A_{{#1},{#2}}$}}

\def \blue[#1]{\textcolor{blue}{#1}}

\title[Eisenbud--Goto conjecture for projectively normal varieties]{The Eisenbud--Goto conjecture for projectively normal varieties with mild singularities}
\author{Jong In Han}
\address{Jong In Han, School of Mathematics, Korea Institute for Advanced Study (KIAS), 85 Hoegi-ro, Dongdaemun-gu, Seoul, 02455, Republic of Korea}
\email{jihan09@kias.re.kr}


\begin{document}
\begin{abstract}
	For a nondegenerate projective variety $X$, the Eisenbud--Goto conjecture asserts that $\operatorname{reg}X\leq\operatorname{deg}X-\operatorname{codim}X+1$. Despite the existence of counterexamples, identifying the classes of varieties for which the conjecture holds remains a major open problem. In this paper, we prove that the Eisenbud--Goto conjecture holds for $2$-very ample projectively normal varieties with factorial, rational, hypersurface singularities and isolated Gorenstein singularities.
\end{abstract}
\maketitle

\section{Introduction}

In this paper, we work over an algebraically closed field $\k$ of characteristic zero and require varieties to be integral.
Let $X\subseteq\p^r$ be a nondegenerate projective variety.
For a coherent sheaf $\f$ on $\p^r$, we say \emph{$\f$ is $m$-regular} if 
\[
H^i(\p^r,\f(m-i))=0
\]
for all $i\geq 1$.
The regularity $\reg \f$ of $\f$ is defined as the least integer $m$ such that $\f$ is $m$-regular.
The regularity $\reg X$ of $X$ is defined as the regularity of its ideal sheaf $\I_X$.
Regularity is one of the invariants measuring the complexity of projective varieties.
In 1984, D. Eisenbud and S. Goto conjectured a bound on the regularity of projective varieties.
\begin{conjecture}[The Eisenbud--Goto conjecture, {\cite{MR741934}}]
	Let $X\subseteq\p^r$ be a nondegenerate projective variety.
	Then
	\begin{equation}\label{EG_bound}
		\reg X\leq \deg X-\codim X+1.
	\end{equation}
\end{conjecture}
It holds for arithmetically Cohen--Macaulay varieties.
In 1983, L. Gruson, R. Lazarsfeld, and C. Peskine \cite{MR704401} showed the conjecture for curves.
For the case of smooth surfaces, H. C. Pinkham \cite{MR818356} and R. Lazarsfeld \cite{MR894589} proved the conjecture.
In 1990, Z. Ran \cite{MR1064868} asserted that the conjecture holds for most smooth threefolds of codimension at least 6 using local differential geometric methods.
In 1998, I. Peeva and B. Sturmfels \cite{MR1649322} proved the conjecture for projective toric varieties of codimension two.
In 1998, S. Kwak \cite{MR1620706} proved the weak bound
\begin{equation}\label{weaker_bound}
	\reg X\leq \deg X-\codim X+2
\end{equation}
for smooth threefolds and $\reg X\leq \deg X-\codim X+5$ for smooth fourfolds.
In 1999, S. Kwak \cite{MR1679165} showed \eqref{EG_bound} for smooth threefolds in $\p^5$.
For smooth projective varieties of dimension at most 14, L. Chiantini, N. Chiarli, and S. Greco \cite{MR1783985} showed
\[
	\reg X\leq \deg X-\codim X+1+\frac{(n-2)(n-1)}{2},
\]
hence improving the bound for smooth fourfolds to $\reg X\leq \deg X-\codim X+4$.
Here, the condition $\dim X\leq 14$ comes from Mather's theorem \cite{MR362393}.
In 2015, W. Niu \cite{MR3369342} showed \eqref{EG_bound} for normal surfaces with rational, Gorenstein elliptic, log canonical singularities.
In 2022, W. Niu and J. Park \cite{MR4395951} showed the weak bound \eqref{weaker_bound} for normal threefolds with rational singularities.

In 2018, J. McCullough and I. Peeva \cite{MR3758150} published a surprising result that there exist counterexamples to the Eisenbud--Goto conjecture.
Their construction uses Rees-like algebras and step-by-step homogenization.
In particular, they constructed a singular threefold of degree 31 and regularity 38 in $\p^5$.
Furthermore, they showed there does not exist a polynomial function of $\deg X$ bounding $\reg X$.
In 2024, the author and S. Kwak \cite{MR4771231} constructed surface counterexamples using binomial rational maps between projective spaces.
In 2025, J. Choe \cite{MR4859414} constructed counterexamples using the unprojection method.

Although counterexamples exist, it is still a very interesting problem to identify the classes of varieties for which the conjecture holds.
As J. McCullough and I. Peeva stated in \cite{MR4195747}, the following are interesting cases where the conjecture may hold:
\begin{enumerate}
	\item smooth projective varieties;
	\item projectively normal varieties;
	\item projective toric varieties; and
	\item projective singular surfaces.
\end{enumerate}
Case (4) was ruled out in \cite{MR4771231}.

On the other hand, the validity of the Eisenbud--Goto conjecture for $X$ is equivalent to the following two conditions:
\begin{enumerate}
	\item for all $m\geq \deg X-\codim X$, the variety $X\subseteq\p^r$ is $m$-normal, i.e., $H^0(\p^r,\o_{\p^r}(m))\to H^0(X,\o_X(m))$ is surjective; and
	\item $\reg \o_X\leq \deg X-\codim X$.	
\end{enumerate}
In 2014, A. Noma \cite{MR3217694} showed
\begin{equation}\label{ox_conj}
	\reg \o_X\leq \deg X-\codim X
\end{equation}
for smooth projective varieties except when $X$ is a scroll over a smooth projective curve.
Here, a \emph{scroll over a smooth curve $C$} is defined as the image of $\p_C(\mathcal{E})$ for some locally free sheaf $\mathcal{E}$ over $C$ under a birational morphism defined by a sublinear system of $\lvert\o_{\p_C(\mathcal{E})}(1)\rvert$.
In 2020, S. Kwak and J. Park \cite{MR4057491} showed \eqref{ox_conj} for smooth scrolls over smooth projective curves.
Therefore \eqref{ox_conj} holds for smooth projective varieties.
This implies that the Eisenbud--Goto conjecture holds for smooth projectively normal varieties.
In 2021, J. Moraga, J. Park, and L. Song \cite{MR4289907} showed \eqref{ox_conj} for normal projective varieties with isolated $\q$-Gorenstein singularities.
Hence the Eisenbud--Goto conjecture holds for projectively normal varieties with isolated $\q$-Gorenstein singularities.

In this paper, we prove that the Eisenbud--Goto conjecture holds for 2-very ample projectively normal varieties with mild singularities.
Regarding mild singularities, we allow two types of singularities: one is the factorial, rational, hypersurface singularity and the other is the isolated Gorenstein singularity.
Specifically, we prove the following.

\begin{thmalphabetintro}\label{intro_EG_for_PN}
	Let $X\subseteq\p^r$ be a nondegenerate $2$-very ample projectively normal variety with factorial, rational, hypersurface singularities and isolated Gorenstein singularities.
	Then
	\[
		\reg X\leq \deg X-\codim X+1.
	\]
\end{thmalphabetintro}

This paper is organized as follows.
In \Cref{sec: preliminaries}, we review some known results.
In \Cref{sec: Kodaira}, we discuss vanishing theorems for singular varieties.
In \Cref{sec: birational double point formula}, we derive the birational double point formula for factorial points.
In \Cref{sec: the main theorem}, we prove the main theorem.

\section*{Acknowledgement}
The author was supported by a KIAS Individual Grant (MG101401) at Korea Institute for Advanced Study.
The author thanks Professor Sijong Kwak for his encouragement.

\section{Preliminaries}\label{sec: preliminaries}
In this section, we collect some known results for later discussion.
Let $X\subseteq\p^r$ be a nondegenerate projective variety of dimension $n$ and codimension $e$.
We recall the definition of $k$-very ampleness, as the main theorem relies on this positivity condition.

\begin{definition}
	Let $\l$ be a line bundle on $X$.
	We say \emph{$\l$ is $k$-very ample} if
	\[
		H^0(X,\l)\to H^0(\xi,\l|_\xi)
	\]
	is surjective for any subscheme $\xi$ of dimension 0 and length $k+1$.
	When $X\subseteq\p^r$ is embedded by the complete linear system $\lvert \l\rvert$, we say \emph{$X$ is $k$-very ample} if $\l$ is $k$-very ample.
\end{definition}

On locally factorial varieties, reflexive sheaves of rank 1 are invertible.
This fact is used in the proof of \Cref{singular_birational_double_point_formula}.
\begin{proposition}[{\cite[Proposition 1.9]{MR597077}}]\label{reflexive_on_fac}
	On integral Noetherian locally factorial schemes, any reflexive sheaf of rank 1 is invertible.
\end{proposition}

The Zariski--Fujita theorem states that if the base locus of a line bundle is finite, then the line bundle is semiample.

\begin{theorem}[{\cite{MR141668}, also cf. \cite[Remark 2.1.31]{MR2095471}}]\label{ZFT}
	Let $\l$ be a line bundle on $X$ with
	\[
		\dim\bs\lvert\l\rvert\leq 0.
	\]
	Then $\l$ is semiample.
\end{theorem}

Now we review the results by A. Noma.
First, we denote the locus of smooth points of $X$ that induce a nonbirational inner projection to $\p^{r-1}$ by
\[
\C(X):=\{x\in X_{\sm}\mid \leng(X\cap \la x, y\ra)\geq 3\text{ for general }y\in X\}
\]
where $\la \cdot\ra$ denotes the linear span.

It is worth noting that the closure of this locus is a union of linear subspaces.

\begin{theorem}[{\cite[Theorem 5]{MR2645037}}]\label{C(X)}
	Let $X\subseteq\p^r$ be a nondegenerate projective variety of codimension $e\geq 2$ and $\Lambda$ be an irreducible component of $\C(X)$. Then $\overline{\Lambda}$ is linear and
	\[
	\dim \C(X)\leq\min\{n-1, \dim X_\sing+2\}.
	\]
\end{theorem}

Let $x_1,\cdots,x_m$ be general points in $X$ and $Z:=\{x_1,\cdots,x_m\}$ where $1\leq m\leq e-1$.
Denote by $\pi:X\setminus \la Z\ra\to\p^{r-m}$ the linear projection of $X\subseteq\p^r$ from $\la Z\ra$ and the closure of the image of $\pi$ by $X_Z$.
The following theorem shows the map $\pi':X\setminus \la Z\ra\to X_Z$ called the \emph{induced projection} is isomorphic at $x$ if $x\in X_\sm\setminus\C(X)$.
Note that we pick general points $x_1,\cdots,x_m$ when $x$ is given, so one can not expect that $\pi'$ is isomorphic on the entire $X_\sm\setminus\C(X)$.
\begin{theorem}[{\cite[Lemma 2.2]{MR3217694}}]\label{isom_at_smooth}
	Assume as above and $e\geq 2$.
	Suppose $x\in X_\sm\setminus\C(X)$.
	Then $\pi':X\setminus \la Z\ra\to X_Z$ is isomorphic at $x$.
\end{theorem}

Denote by $E_{x_1,\cdots,x_m}(X)$ the union of positive dimensional fibers of $\pi$, i.e.,
\[
E_{x_1,\cdots,x_m}(X):=\{x\in X\setminus \la Z\ra\mid \dim\pi^{-1}(\pi(x))\geq 1\}.
\]

\begin{definition}
	We say \emph{$X$ satisfies $(E_m)$} if
	\[
	\dim E_{x_1,\cdots,x_m}(X)\geq n-1,
	\]
	or equivalently, $\dim E_{x_1,\cdots,x_m}(X)=n-1$, where $x_1,\cdots,x_m\in X$ are general points ($1\leq m\leq e-1$).
\end{definition}

Let $\bl_ZX$ be the blowing up of $X$ along $Z$.
Then we get the map $\tilde{\pi}:\bl_ZX\to\p^{r-m}$ called the \emph{extended projection}.
With this notation, the following lemma holds.

\begin{lemma}[{\cite[Lemma 1.3]{MR3217694}}]\label{extended_proj_lemma}
	Assume as above and $n\geq 2$.
	Then
	\[
	\tilde{E}_{x_1,\cdots,x_m}(X):=\{x\in \bl_ZX\mid \dim\tilde{\pi}^{-1}(\tilde{\pi}(x))\geq 1\}
	\]
	is closed in $\bl_ZX$ and $\dim E_{x_1,\cdots,x_m}(X)\geq n-1$ if and only if $\dim \tilde{E}_{x_1,\cdots,x_m}(X)\geq n-1$.
\end{lemma}

Using the condition $(E_m)$, Noma established the positivity of the double point divisor for inner projections.

\begin{theorem}[{\cite[Theorem 1]{MR3217694}}]\label{positivity_dpd}
	Let $X\subseteq\p^r$ be a nondegenerate smooth projective variety of dimension $n\geq 2$, codimension $e\geq 2$, and degree $d$ not satisfying $(E_m)$ where $1\leq m\leq e-1$.
	Then
	\[
	\bs\lvert \o_X(d-n-m-2)\otimes \omega_X^{-1}\rvert\subseteq \C(X)
	\]
	and $\o_X(d-n-m-1)\otimes \omega_X^{-1}$ is very ample on $X\setminus\C(X)$.
	Furthermore if $\dim \C(X)\leq 0$, then $\o_X(d-n-m-2)\otimes \omega_X^{-1}$ is semiample so that $\o_X(d-n-m-1)\otimes \omega_X^{-1}$ is ample.
\end{theorem}

We define a scroll over a smooth curve as follows.

\begin{definition}
	A \emph{scroll over a smooth curve $C$} is the image of $\p_C(\mathcal{E})$ for some locally free sheaf $\mathcal{E}$ over $C$ under a birational morphism defined by a sublinear system of $\lvert\o_{\p_C(\mathcal{E})}(1)\rvert$.
\end{definition}

Noma classified projective varieties satisfying $(E_m)$ for $1\leq m\leq e-1$.
\begin{theorem}[{\cite[Theorem 3]{MR3217694}}]\label{classification_Em}
	Let $X\subseteq\p^r$ be a nondegenerate projective variety of dimension $n\geq 2$ and codimension $e$.
	\begin{enumerate}
		\item When $e\geq 2$ and $X$ satisfies $(E_1)$, $X$ is projectively equivalent to a scroll over a smooth curve.
		If furthermore $X$ is smooth, $X$ is projectively equivalent to $\p_C(\mathcal{E})\subseteq\p^r$ where $\mathcal{E}$ is a locally free sheaf on a smooth projective curve.
		\item When $e\geq 3$ and $X$ satisfies $(E_2)$ but not $(E_1)$, $X$ is projectively equivalent to a cone over the Veronese surface $\nu_2(\p^2)\subseteq\p^5$.
		If furthermore $X$ is smooth, $X$ is projectively equivalent to the Veronese surface $\nu_2(\p^2)\subseteq\p^5$.
		\item When $e\geq 4$, $X$ satisfies $(E_m)$ for some $3\leq m\leq e-1$ if and only if $X$ satisfies $(E_1)$.
	\end{enumerate}
\end{theorem}

In particular, if $n\geq 2$ and $X$ is not projectively equivalent to
\begin{enumerate}
	\item a scroll over a smooth curve $C$; and
	\item a cone over the Veronese surface $\nu_2(\p^2)\subseteq\p^5$,
\end{enumerate}
then $X$ does not satisfy $(E_m)$ for any $1\leq m\leq e-1$.

\section{Vanishing theorems for singular varieties}\label{sec: Kodaira}
Let $X\subseteq\p^r$ be a nondegenerate projective variety of dimension $n$. 
We denote the normalized dualizing complex of $X$ by $\omega_X^\bullet$ (cf. \cite{MR222093}) and define the dualizing sheaf as $\omega_X:=\mathcal{H}^{-n}(\omega_X^\bullet)$, the lowest cohomology of $\omega_X^\bullet$.
When $X$ is Cohen--Macaulay, this is the only cohomology that does not vanish.
D. Arapura, D. B. Jaffe, and L. Song \cite{MR952313,AS18} studied Kodaira vanishing theorem for singular projective varieties.
One of their vanishing results for singular projective varieties is the following.

\begin{theorem}[{\cite[Lemma 3.3]{AS18}}]\label{vanishing1}
	Let $X\subseteq\p^r$ be a projective variety and $\l$ be a nef and big line bundle on $X$.
	Then
	\[
		H^i(X,\H^j(\omega_X^\bullet)\otimes \l)=0
	\]
	for any positive integer $i>\dim X_\sing$ and any $j\in \z$.
\end{theorem}

The key part of the proof is to utilize the following vanishing theorem which is a special case of \cite[Th\'{e}or\`{e}me 3.1.A]{MR907924} or a consequence of Grauert--Riemenschneider vanishing and Kawamata--Viehweg vanishing.

\begin{theorem}[{cf. \cite[Th\'{e}or\`{e}me 3.1.A]{MR907924}}]\label{vanishing2}
	Let $X\subseteq\p^r$ be a projective variety, $\l$ a nef and big line bundle on $X$, and $f:Y\to X$ be a resolution of singularities. Then
	\[
		H^i(X,f_*\omega_Y\otimes \l)=0
	\]
	for all $i\geq 1$.
\end{theorem}

We define the locus
\[
X_\cm := \{x\in X\mid \text{$X$ is Cohen--Macaulay at $x$}\}.
\]
This locus is open in $X$ (cf. \cite[Exercise 24.2]{MR1011461}).
We also define the locus
\[
X_\rat := \{x\in X\mid \text{$X$ has at most a rational singularity at $x$}\}.
\]
This locus is open in $X$ \cite{MR501926}.

With these loci, one can express the vanishing theorem for singular varieties (\Cref{vanishing1}) as follows.
\Cref{vanishing_rat} is used in the proof of \Cref{EG_bound_for_PN}.

\begin{lemma}\label{vanishing_rat}
	Let $X\subseteq\p^r$ be a projective variety of dimension $n$ and $\l$ be a nef and big line bundle on $X$.
	Then
	\[
	H^i(X,\mathcal{H}^j(\omega_X^\bullet)\otimes \l)=0
	\]
	for any nonnegative integer $i> \dim (X\setminus X_\cm)$ when $j\geq -n+1$ and for any positive integer $i> \dim (X\setminus X_\rat)$ when $j= -n$.
	In particular, it holds for any positive integer $i> \dim (X\setminus X_\rat)$ and any $j\in \z$.
\end{lemma}
\begin{proof}
	Denote $n:=\dim X$.
	For $j\geq -n+1$, the sheaf $\mathcal{H}^j(\omega_X^\bullet)$ vanishes on $X_\cm$.
	Hence it is supported on $X\setminus X_\cm$ so that
	\[
		H^i(X,\mathcal{H}^j(\omega_X^\bullet)\otimes \l)=0
	\]
	for any nonnegative integer $i>\dim(X\setminus X_\cm)$ and any integer $j\geq -n+1$.
	
	Now let $j=-n$.
	For a resolution of singularities $f:Y\to X$, we have the injective trace map $\operatorname{Trace}_{Y/X}:f_*\omega_Y\to \omega_X$ (cf. \cite[Proposition 5.77]{MR1658959}).
	Also,
	\[
		H^i(X,f_*\omega_Y\otimes \l)=0
	\]
	for all $i>0$ (\Cref{vanishing2}).
	On the other hand, the cokernel of $\operatorname{Trace}_{Y/X}$ is supported on $X\setminus X_\rat$.
	Hence $H^i(X,\omega_X\otimes\l)=0$ for all positive $i>\dim(X\setminus X_\rat)$.
	The last statement follows from $X\setminus X_\cm \subseteq X\setminus X_\rat$.
\end{proof}

\section{The birational double point formula}\label{sec: birational double point formula}

The purpose of this section is to prove \Cref{singular_birational_double_point_formula}, which is used in the proof of \Cref{EG_bound_for_PN}.
R. Lazarsfeld \cite{MR2095472} and A. Noma \cite{MR3217694} stated the birational double point formula which is derived from the double point formula \cite[Theorem 9.3]{MR1644323}.

\begin{lemma}[{\cite[Lemma 10.2.8]{MR2095472} and \cite[Lemma 2.1]{MR3217694}}]\label{smooth_birational_double_point_formula}
	Let $f:X\to P$ be a morphism between smooth projective varieties where $\dim X=n$ and $\dim P=n+1$.
	Suppose $f$ factors through a surjective birational morphism $g:X \twoheadrightarrow Y$ for a subvariety $Y\subseteq P$ of codimension 1.
	Then for some effective divisors $D$ and $E$ on $X$, it holds that
	\[
	D-E\equiv -K_X+f^*(K_P+Y)
	\]
	and $E$ is exceptional for $g$.
	Furthermore, if $f$ is isomorphic at a point $x\in X$, then $x\notin \supp(D-E)$.
\end{lemma}

In \cite{MR4289907}, the authors derived the birational double point formula as follows.
\begin{lemma}[{\cite[Lemma 2.3]{MR4289907}}]
	Let $f:X\to \p^{n+1}$ be a morphism from a normal projective variety of dimension $n$.
	Suppose $f$ factors through a surjective birational morphism $g:X \twoheadrightarrow Y$ for a hypersurface $Y\subseteq \p^{n+1}$.
	Then for some effective Weil divisors $D$ and $E$ on $X$, it holds that
	\[
	D-E\equiv -K_X+f^*(K_P+Y)
	\]
	and $E$ is exceptional for $g$.
	Furthermore, if $f$ is isomorphic at a smooth point $x\in X_\sm$, then $x\notin \supp(D-E)$.
\end{lemma}

Here, \emph{$E$ is exceptional for $g$} means that the image of every irreducible component of $E$ has dimension lower than $E$.
The divisor $D$ is called the \emph{double point divisor}.
The positivity of the double point divisor has been investigated in several works \cite{MR1253986,MR1422899,MR3217694,CP24}.

In the following, we derive the birational double point formula for factorial points to prove \Cref{EG_bound_for_PN}.
To state it, we define the locus
\[
X_\fac := \{x\in X\mid \text{$X$ is factorial at $x$}\}.
\]
This locus is open in $X$ \cite{MR749675}, also cf. \cite{BGS11}.

\begin{lemma}\label{singular_birational_double_point_formula}
	Let $f:X\to Y$ be a birational morphism between Gorenstein projective varieties of dimension $n$ not admitting any exceptional divisor. Suppose $X$ is normal. If $f$ is isomorphic at a factorial point $x\in X_{\fac}$, then $\lvert \omega_X^{-1}\otimes f^*\omega_Y\rvert\neq \varnothing$ and $x\notin \bs \lvert\omega_X^{-1}\otimes f^*\omega_Y\rvert$.
\end{lemma}

\begin{proof}
	Let $i:X_{\fac}\to X$ be the immersion.
	Denote $g:=f\circ i$ so that we have the following commutative diagram.
	
	\[
	\begin{tikzcd}
		X_{\fac}\arrow[d,"i"]\arrow[rd,"g"]&\\
		X\arrow[r,"f"]&Y
	\end{tikzcd}
	\]
	
	Note that $g$ is birational and does not admit an exceptional divisor.
	By \Cref{reflexive_on_fac}, the sheaf $(g^*g_*\omega_{X_{\fac}})^{\vee\vee}$ is a line bundle since it is reflexive of rank 1 and $X_\fac$ is locally factorial.
	The counit map $g^*g_*\omega_{X_{\fac}}\to \omega_{X_{\fac}}$ factors through $(g^*g_*\omega_{X_{\fac}})^{\vee\vee}$.
	The map $(g^*g_*\omega_{X_{\fac}})^{\vee\vee}\to \omega_{X_{\fac}}$ is injective since $g$ is birational and $(g^*g_*\omega_{X_{\fac}})^{\vee\vee}$ is torsion free.
	Hence $(g^*g_*\omega_{X_{\fac}})^{\vee\vee}=\omega_{X_{\fac}}(-E)$ for some effective divisor $E$ that is exceptional for $g$.
	As $g$ does not admit an exceptional divisor, we get
	\[
		(g^*g_*\omega_{X_{\fac}})^{\vee\vee}=\omega_{X_{\fac}}.
	\]
	
	On the other hand, as $f$ is a birational morphism between projective varieties, there exists the trace map $f_*\omega_X\to\omega_Y$ that is injective.
	Note that $i_*\omega_{X_{\fac}}=i_*i^*\omega_X=\omega_X$ as $X$ is normal and $\dim (X\setminus X_{\fac})\leq n-2$.
	Hence
	\[
		g_*\omega_{X_\fac}\otimes \omega_Y^{-1}=f_*i_*\omega_{X_{\fac}}\otimes \omega_Y^{-1}=f_*\omega_X\otimes \omega_Y^{-1}\hookrightarrow \o_Y.
	\]
	This induces the map
	{\small
	\[
		\omega_{X_{\fac}}\otimes g^*\omega_Y^{-1}=(g^*g_*\omega_{X_\fac})^{\vee\vee}\otimes g^*\omega_Y^{-1}=(g^*g_*\omega_{X_\fac}\otimes g^*\omega_Y^{-1})^{\vee\vee}\to (g^*\o_Y)^{\vee\vee}=\o_{X_\fac}
	\]
	}
	which is injective as $g$ is birational.
	Therefore $\omega_{X_\fac}\otimes g^*\omega_Y^{-1}$ is an invertible ideal sheaf on $X_\fac$ so that
	\[
		\omega_{X_\fac}\otimes g^*\omega_Y^{-1}=\o_{X_\fac}(-D_0)
	\]
	for some effective divisor $D_0$.
	
	As $X$ is normal, we have $\dim (X\setminus X_\fac)\leq n-2$ so that the divisor $D_0$ extends to an effective Weil divisor $D$ on $X$. Also note that $i_*(\omega_{X_\fac}^{-1})=i_*i^*(\omega_{X}^{-1})=\omega_X^{-1}$.
	Applying $i_*$ to $\o_{X_\fac}(D_0)=\omega_{X_\fac}^{-1}\otimes g^*\omega_Y$, we get
	\[
		\o_X(D)=i_*(\omega_{X_\fac}^{-1}\otimes i^*f^*\omega_Y)=i_*(\omega_{X_\fac}^{-1})\otimes f^*\omega_Y=\omega_X^{-1}\otimes f^*\omega_Y.
	\]
	Hence $D\in \lvert\omega_X^{-1}\otimes f^*\omega_Y\rvert$.
	
	Suppose $f$ is isomorphic at a factorial point $x\in X_\fac$. Then the above injective map
	\[
		\omega_{X_\fac}\otimes g^* \omega_Y^{-1}\to\o_{X_\fac}
	\]
	is isomorphic at $x$ so that $x\notin D_0$.
	As $x\in X_\fac$, this implies $x\notin D$. Therefore $x\notin \bs\lvert\omega_X^{-1}\otimes f^*\omega_Y\rvert$.
\end{proof}

\section{The Eisenbud--Goto regularity bound}\label{sec: the main theorem}

In \cite{MR3217694}, Noma bounded the regularity of $\o_X$ using his results on the positivity of the double point divisor (\Cref{positivity_dpd}), the classification of projective varieties satisfying $(E_m)$ (\Cref{classification_Em}), and the characterization of smooth projective varieties with $\dim\C(X)\geq 1$ (\cite[Corollary 6.2]{MR3748569}).

\begin{theorem}[{\cite[Corollary 5]{MR3217694}}]
	Let $X\subseteq\p^r$ be a nondegenerate smooth projective variety of dimension $n\geq 2$, codimension $e$, and degree $d$ that is not projectively equivalent to a scroll over a smooth projective curve.
	Then
	\[
		\reg\o_X\leq d-e.
	\]
\end{theorem}

The remaining case was solved by S. Kwak and J. Park \cite{MR4057491}.

\begin{theorem}[{\cite[Proposition 3.6]{MR4057491}}]\label{scroll_reg}
	Let $X\subseteq\p^r$ be a nondegenerate smooth scroll over a smooth projective curve $C$ of dimension $n\geq 2$, codimension $e$, and degree $d$.
	Then
	\begin{enumerate}
		\item $\reg\o_X\leq d-e$;
		\item $\reg\o_X=d-e$ if and only if $X$ is a rational normal scroll; and
		\item $\reg\o_X=d-e-1$ if and only if $X$ is an isomorphic projection of rational normal scroll from a point or an elliptic normal surface scroll with $d=e+3$.
	\end{enumerate}
\end{theorem}

Furthermore, J. Moraga, J. Park, and L. Song \cite{MR4289907} derived the result for normal projective varieties with isolated $\q$-Gorenstein singularities.

\begin{theorem}[{\cite[Theorem 1]{MR4289907}}]\label{normal_isol}
	Let $X\subseteq\p^r$ be a nondegenerate normal variety with isolated $\q$-Gorenstein singularities of codimension $e$ and degree $d$.
	Then $$\reg\o_X\leq d-e.$$
\end{theorem}

To prove the main theorem, we introduce two lemmas: \Cref{isom_onestep} and \Cref{isom_inductive}.
These lemmas follow the argument of \cite[Lemma 2.2]{MR3217694} (\Cref{isom_at_smooth}) and extend the statement to singular points $x$ with $\dim T_xX\leq n+1$.
\begin{lemma}\label{isom_onestep}
	Let $X\subseteq\p^r$ be a nondegenerate projective variety of dimension $n$ and codimension $e\geq 2$.
	For a point $x\in X$, suppose that $\dim T_xX\leq n+1$ and $\leng(X\cap \la x,x_1\ra)=2$ for the general point $x_1\in X$.
	Let $X_{x_1}$ be the closure of the image of the linear projection $\pi:X\setminus \{x_1\}\to \p^{r-1}$ of $X$ from $x_1$.
	Then the induced projection $\pi':X\setminus \{x_1\}\to X_{x_1}$ is isomorphic at $x$.
	Furthermore, if $e\geq 3$, then $\leng (X_{x_1}\cap \la \pi(x),y_\gen\ra)=2$ for the general point $y_\gen\in X_{x_1}$.
\end{lemma}
\begin{proof}
	Note that $x$ is not a vertex of $X$ by the assumption.
	Let $b:\bl_{x_1}X\to X$ be the blowing up of $X$ along $x_1$ with the exceptional divisor $E_1$ and $\tilde{x}$ be the preimage of $x$ via $b$.
	By the generic smoothness theorem (cf. \cite[Lemma III.10.5]{MR0463157}), the inner projection of $X$ from $x$ is unramified at $x_1$, so we have $x\notin T_{x_1}X$.
	Thus the extended projection $\bl_{x_1}X\to\p^{r-1}$ separates $\tilde{x}$ from any point in $E_1$.
	Also as $\leng (X\cap\la x,x_1\ra)=2$, the extended projection $\bl_{x_1}X\to\p^{r-1}$ separates $\tilde{x}$ from any other point in $\bl_{x_1}X\setminus E_1$.
	Therefore $\bl_{x_1}X\to\p^{r-1}$ separates $\tilde{x}$ from any other point in $\bl_{x_1}X$.
	Furthermore, $x_1\notin T_xX$ by the assumption that $\dim T_xX\leq n+1$ and $e\geq 2$.
	Therefore $\bl_{x_1}X\to\p^{r-1}$ is an embedding near $\tilde{x}$, so $\pi':X\setminus \{x_1\}\to X_{x_1}$ is isomorphic at $x$.
	
	Now we assume $e\geq 3$ to show the last statement.
	It is enough to show $\leng(X_{x_1}\cap \la \pi(x), \pi(y)\ra)=2$ for the general point $y\in X$.
	We suppose the contrary, i.e.,
	\[
	\leng(X_{x_1}\cap \la \pi(x), \pi(y)\ra)\geq 3
	\]
	for the general point $y\in X$ or $\pi(x)$ is a vertex point of $X_{x_1}$.
	The latter case cannot occur since $\dim T_{\pi(x)}X_{x_1}\leq n+1$ as $\pi':X\setminus \{x_1\}\to X_{x_1}$ is isomorphic at $x$.
	
	We claim that there exists a point $z\in X\setminus \{x_1\}$ such that $\pi(z)\in \la \pi(x),\pi(y)\ra$, $\pi(z)\neq \pi(x)$, and $\pi(z)\neq \pi(y)$.
	As $\dim T_{\pi(x)}X_{x_1}\leq n+1$, we have $\dim \la T_{\pi(x)}X_{x_1}, x_1\ra\leq n+2$.
	As $e\geq 3$, this implies $y\notin \la T_{\pi(x)}X_{x_1}, x_1\ra$ so that $\pi(y)\notin T_{\pi(x)}X_{x_1}$.
	Also, by the generic smoothness theorem, the inner projection of $X_{x_1}$ from $\pi(x)$ is unramified at $\pi(y)$, so we have $\pi(x)\notin T_{\pi(y)}X_{x_1}$.
	Furthermore, as $\dim \la T_{x_1}X,x\ra=n+1$, we have $y\notin \la T_{x_1}X,x\ra$ so that any points in $E_1$ do not map to $\la \pi(x),\pi(y)\ra$.
	Therefore the claim holds.
	
	By the claim, the points $x_1,x,y,z$ are distinct points on a $2$-plane.
	Furthermore, $z\notin \la x,x_1\ra$ and $z\notin \la x,y\ra$ as $$\leng (X\cap \la x,x_1\ra)=\leng (X\cap \la x,y\ra)=2$$ by the assumption.
	Let $X_x$ be the closure of the image of the inner projection $\pi_x:X\setminus\{x\}\to\p^{r-1}$ of $X$ from $x$.
	Then $\la \pi_x(x_1),\pi_x(y)\ra$ is the general secant line to $X_x$, but is also a trisecant line as it contains $\pi_x(z)$.
	This contradicts the general position lemma (cf. \cite[Lemma 1.2]{MR3217694}).
\end{proof}

\begin{lemma}\label{isom_inductive}
	Let $X\subseteq\p^r$ be a nondegenerate projective variety of dimension $n$ and codimension $e\geq 2$.
	For a point $x\in X$, suppose that $\dim T_xX\leq n+1$ and $\leng(X\cap \la x,x_\gen\ra)=2$ for the general point $x_\gen\in X$.
	Let $x_1,\cdots,x_{e-1}\in X$ be general points and $Z:=\{x_1,\cdots,x_{e-1}\}$.
	Denote by $X_Z$ the closure of the image of the linear projection $\pi:X\setminus \la Z \ra\to \p^{n+1}$ of $X$ from $\la Z \ra$.
	Then the induced projection $\pi':X\setminus \la Z \ra\to X_Z$ is isomorphic at $x$.
\end{lemma}
\begin{proof}
	This follows by applying \Cref{isom_onestep} inductively.
	Note that when $e\geq 3$, \Cref{isom_onestep} ensures that the condition
	\begin{equation}\label{analogous_CX_condition}
		\leng(X\cap \la x,x_\gen\ra)=2\text{ for the general point }x_\gen\in X
	\end{equation}
	holds also in the next step.
	When $e=2$, \Cref{isom_onestep} does not ensure that \eqref{analogous_CX_condition} holds for the image, but it holds for $X$ itself, which suffices to apply the projection.
\end{proof}

In the following, the proof follows the argument of Noma (\Cref{positivity_dpd}) which considers smooth projective varieties.
We extend it to the singular case under consideration.

\begin{theorem}\label{EG_bound_for_PN}
	Let $X\subseteq\p^r$ be a nondegenerate $2$-very ample normal Gorenstein variety.
	Suppose that
	\begin{enumerate}
		\item $X$ has factorial rational singularities; and
		\item $\dim T_pX\leq \dim X+1$ for $p\in X$
	\end{enumerate}
	possibly except some finite points.
	Then
	\[
	\reg \o_X\leq \deg X-\codim X.
	\]
\end{theorem}
\begin{proof}
	Denote $n:=\dim X$, $e:=\codim X$, and $d:=\deg X$.
	We assume $n\geq 2$ and $e\geq 2$ since the case $n=1$ is well known \cite{MR704401} and the case $e=1$ is trivial.
	By the assumption, $X$ is neither a cone nor a scroll over a smooth curve.
	If $X$ is the Veronese surface $\nu_2(\p^2)$, then the result is well known as it has minimal degree.
	Thus we may assume $X$ is not $\nu_2(\p^2)$.
	Then $X$ does not satisfy ($E_{e-1}$) by \Cref{classification_Em}.
	
	Let $Y$ be a cofinite subset of $X_\sing$ such that
	\begin{enumerate}
		\item $X$ has factorial rational singularities on $Y$; and
		\item $\dim T_pX\leq n+1$ for all $p\in Y$.
	\end{enumerate}
	Such $Y$ must exist by the assumption.
	Pick a point $x\in X_\sm\cup Y$ and $e-1$ general points $x_1,\cdots,x_{e-1}$ in $X$.
	Denote $Z:=\{x_1,\cdots,x_{e-1}\}$.
	Let $b:\bl_ZX\to X$ be the blowing up of $X$ along $Z$ with exceptional divisors $E_1,\cdots,E_{e-1}$ such that each $E_i$ corresponds to $x_i$.
	Note that $b$ is isomorphic at $x$.
	We denote the preimage of $x$ by $\tilde{x}$, i.e., $b(\tilde{x})=x$.
	Let $\pi:X\setminus \la Z\ra\to\p^{n+1}$ be the linear projection of $X$ from the linear span $\la Z\ra$ of $Z$.
	We denote the closure of the image of $\pi:X\setminus \la Z\ra\to\p^{n+1}$ by $X_Z\subseteq\p^{n+1}$ and denote the induced projection by $\pi':X\setminus\la Z\ra\to X_Z$.
	Also, let $\tilde{\pi}:\bl_ZX\to \p^{n+1}$ be the extended projection and $\tilde{\pi}':\bl_ZX\to X_Z$ be the induced morphism.
	\[
	\begin{tikzcd}
		\bl_ZX\arrow[d,"b"]\arrow[r,"\tilde{\pi}'"]\arrow[rr,"\tilde{\pi}",bend left=24]& X_Z\arrow[r,hook] &\p^{n+1}\\
		X\arrow[ru,"\pi'",dashed]\arrow[r,hook]\arrow[rru,"\pi",dashed]& \p^r\arrow[ru,dashed]&
	\end{tikzcd}
	\]
	
	The morphism $\tilde{\pi}'$ is birational as $\pi'$ is.
	The birational map $\pi'$ does not admit an exceptional divisor as $X$ does not satisfy $(E_{e-1})$.
	Then $\tilde{\pi}'$ also does not admit an exceptional divisor by \Cref{extended_proj_lemma}.
	
	The projection $\pi':X\setminus\la Z\ra\to X_Z$ is isomorphic at $x$ by \Cref{isom_inductive}.
	Also, by the generic smoothness theorem, we have $x\notin T_{x_1}X=\la T_{x_1}X,x_1\ra$.
	As $\dim \la x,T_{x_1}X,x_1\ra=
	n+1$, we have $x_2\notin \la x,T_{x_1}X,x_1\ra$.
	Thus $x\notin \la T_{x_1}X,x_1,x_2\ra$.
	Repeating this inductively, we get $x\notin \la T_{x_1}X,x_1,\cdots,x_{e-1}\ra$.
	This shows that $x\notin \la T_{x_i}X,x_1,\cdots,x_{e-1}\ra$ for all $i=1,\cdots,e-1$.
	Hence $\tilde{\pi}':\bl_ZX\to X_Z$ separates $\tilde{x}$ from any point in $E:=E_1+\cdots+E_{e-1}$.
	Together with the fact that $\pi':X\setminus\la Z\ra\to X_Z$ is isomorphic at $x$, this implies that $\tilde{\pi}':\bl_ZX\to X_Z$ is isomorphic at $\tilde{x}$.
	
	On the other hand, as $X_Z$ is a hypersurface of $\p^{n+1}$, it is Gorenstein.
	Also, $\bl_Z X$ is a normal Gorenstein variety as $X$ is a normal Gorenstein variety and $x_1,\cdots,x_{e-1}$ are smooth points.
	As $x\in X_\sm\cup Y$, it holds that $X$ is factorial at $x$ so that $\bl_ZX$ is factorial at $\tilde{x}$.
	Hence by \Cref{singular_birational_double_point_formula}, $\tilde{x}\notin \bs\lvert \omega_{\bl_ZX}^{-1}\otimes(\tilde{\pi}')^*\omega_{X_Z}\rvert$.
	Let $\tilde{D}\in \lvert\omega_{\bl_ZX}^{-1}\otimes(\tilde{\pi}')^*\omega_{X_Z}\rvert$ be the divisor not containing $\tilde{x}$.	
	As $b:\bl_ZX\to X$ is isomorphic outside $E:=E_1+\cdots+E_{e-1}$, by taking the isomorphic image of $\tilde{D}|_{\bl_ZX\setminus E}$ and taking the closure, we get an effective divisor $D$ on $X$ not containing $x$ such that
	\[
		\o_X(D)|_{X\setminus Z}=(\omega_X^{-1}\otimes(\pi')^*\omega_{X_Z})|_{X\setminus Z}.
	\]
	This equality extends to the entire $X$ since $X$ is normal with $n\geq 2$, the set $Z$ is finite, and $(\omega_X^{-1}\otimes(\pi')^*\omega_{X_Z})|_{X\setminus Z}$ is a line bundle on $X\setminus Z$.
	Note that $\deg X_Z=d-e+1$.
	Hence
	\[
		\omega_{X_Z}=(\omega_{\p^{n+1}}\otimes \o_{\p^{n+1}}(X_Z))|_{X_Z}=\o_{X_Z}(d-n-e-1),
	\]
	so $D\in \lvert\omega_X^{-1}(d-n-e-1)\rvert$.
	Therefore $x\notin \bs\lvert\omega_X^{-1}(d-n-e-1)\rvert$.
	By repeating this for all $x\in X_\sm\cup Y$, we get
	\[
		\bs\lvert\omega_X^{-1}(d-n-e-1)\rvert\subseteq X_\sing\setminus Y.
	\]
	
	Hence $\omega_X^{-1}(d-n-e-1)$ is semiample by Zariski--Fujita's theorem (\Cref{ZFT}). 
	Therefore $\omega_X^{-1}(d-n-e)$ is ample.
	Note that $\dim(X\setminus X_\rat)\leq 0$ by the assumption.
	Hence by \Cref{vanishing_rat}, we get
	\[
		H^i(X,\o_X(j))=0\text{ for all $1\leq i\leq n$ and $j\geq d-n-e$.}
	\]
	Thus $\o_X$ is $(d-e)$-regular.
\end{proof}

We obtain \Cref{maintext_EG_for_PN} as a corollary.
\begin{thmalphabetmaintext}\label{maintext_EG_for_PN}
	Let $X\subseteq\p^r$ be a nondegenerate $2$-very ample projectively normal variety with factorial, rational, hypersurface singularities and isolated Gorenstein singularities.
	Then
	\[
	\reg X\leq \deg X-\codim X+1.
	\]
\end{thmalphabetmaintext}
\begin{proof}
	As before, denote $e:=\codim X$ and $d:=\deg X$.
	Since any factorial Cohen--Macaulay point is a Gorenstein point (cf. \cite[Corollary 3.3.19]{MR1251956}), $X$ is Gorenstein.
	Hence by \Cref{EG_bound_for_PN}, $\o_X$ is $(d-e)$-regular.
	From the long exact sequence
	{\small
		\[
		\cdots \to H^i(\p^r,\o_{\p^r}(j)) \to H^i(X,\o_X(j)) \to H^{i+1}(\p^r,\I_X(j)) \to H^{i+1}(\p^r,\o_{\p^r}(j)) \to \cdots,
		\]
	}
	it follows that \[H^{i+1}(\p^r,\I_X(d-e-i))=0\text{ for all }i\geq 1.\]
	As $X$ is projectively normal, we have $H^1(\p^r,\I_X(j))=0$ for all $j\geq 0$.
	Therefore, $\I_X$ is $(d-e+1)$-regular.
\end{proof}

\bibliographystyle{amsalpha}
\bibliography{ref.bib}
\end{document}